\documentclass[a4paper,10pt]{article}
\usepackage{amsmath,amscd,amssymb}
\newcommand{\qed}{\hfill$\square$}
\begin{document}
\centerline{\textsc{Espace des modules}}
\centerline{\textsc{des faisceaux de rang 2 semi-stables}}
\centerline{\textsc{de classes de Chern $c_{1}=0$, $c_{2}=2$ et 
$c_{3}=0$}}
\centerline{\textsc{sur la cubique de $\mathbb{P}^{4}$}}
$\ $
\vspace{0.5cm}\\
\centerline{St\'ephane \textsc{Druel}}
$\ $
\vspace{1cm}\\
\textbf{1. Introduction}\\
$\ $
\newline
\indent (1.1) Soient $X\subset\mathbb{P}^{4}$ une hypersurface cubique 
lisse et $\ell\subset X$ une droite de $\mathbb{P}^{4}$. 
Soit $X_{\ell}$ la vari\'et\'e obtenue en \'eclatant $\ell$ dans $X$. 
La projection le long de $\ell$ induit un morphisme 
$X_{\ell}\overset{p}{\longrightarrow}\mathbb{P}^{2}$ dont les fibres 
sont les coniques qui sont coplanaires avec $\ell$. Lorsque la 
droite $\ell$ est g\'en\'erique les fibres de $p$ sont lisses 
ou r\'eunion de deux droites distinctes. Le lieu de 
d\'eg\'en\'erescence de $p$ est alors une courbe 
plane lisse et connexe $C_{0}$ de degr\'e 5. Soit
$C$ la vari\'et\'e des droites contenues dans $X$ et incidentes 
\`a $\ell$. Le morphisme 
$C\longrightarrow C_{0}$ est un rev\^etement \'etale double connexe.
Soit $i$ l'involution \'echangeant les deux feuillets dudit
rev\^etement et notons encore $i$ l'automorphisme 
induit sur la jacobienne $JC$. La vari\'et\'e de Prym associ\'ee 
au rev\^etement $(C,C_{0})$ est alors $P=(Id-i)JC$. C'est une 
vari\'et\'e ab\'elienne principalement polaris\'ee de dimension 5. 
Soient $A_{1}(X)$ le groupe des 1-cycles alg\'ebriques modulo 
l'\'equivalence rationnelle et $A\subset A_{1}(X)$ le sous-groupe 
des cycles alg\'ebriquement \'equivalents \`a z\'ero. 
L'application qui \`a $t\in C$ associe la classe de la droite 
correspondante $z_{t}\subset X_{\ell}$ dans $A$ induit 
un isomorphisme de groupes $P\simeq A$. On d\'emontre que 
pour toute vari\'et\'e lisse $T$ de dimension 
pure $n\ge 1$ et tout $n+1$-cycle $z$ sur $X\times T$ l'application 
d'Abel-Jacobi  
qui \`a $t\in T$ associe la classe du cycle $z_{t}-z_{t_{0}}$ 
dans $P$, o\`u $t_{0}\in T$ est fix\'e, est alg\'ebrique ([Mu]). 
La jacobienne interm\'ediaire de $X$ est d\'efinie par :
$$J(X)=(H^{2,1}(X))^{*}/\alpha(H_{3}(X,\mathbb{Z}))$$
o\`u $\alpha$ est l'application donn\'ee par int\'egration sur les cycles.
C'est une vari\'et\'e ab\'elienne principalement polaris\'ee de dimension 
5 isomorphe \`a la vari\'et\'e de Prym.  
Via cet isomorphisme, l'image du cycle $z_{t}-z_{t_{0}}$ 
par l'application d'Abel-Jacobi est la forme 
lin\'eaire donn\'ee par int\'egration sur le cycle 
$\Gamma$ modulo le groupe $\alpha(H_{3}(X,\mathbb{Z}))$ o\`u 
$\partial\Gamma=z_{t}-z_{t_{0}}$. On 
d\'emontre enfin que   
l'application d'Abel-Jacobi induit un plongement de la surface de Fano 
de $X$ dans $J(X)$.\\
\newline
\indent (1.2) Soient $(X,\mathcal{O}_{X}(1))$ une vari\'et\'e polaris\'ee 
de dimension $n\ge 1$ et $E$ un faisceau coh\'erent sur $X$ de 
rang $r$. La pente $\mu(E)$ 
de $E$ est d\'efinie par la formule :
$$\mu(E)=
\frac{c_{1}(E){c_{1}(\mathcal{O}_{X}(1))}^{n-1}}{r}$$
Le faisceau $E$ est dit 
$\mu$-semi-stable (resp. semi-stable) s'il est sans torsion et si 
pour tout sous-faisceau 
$L\subset E$ de rang $0<r'<r$ on a 
$\mu(L)\le\mu(E)$ 
(resp. $\frac{\chi(L(n))}{r'}\le
\frac{\chi(E(n))}{r}$ pour $n\gg 0$).
Il est dit 
$\mu$-stable (resp. stable) s'il est sans torsion et si 
pour tout sous-faisceau 
$L\subset E$ de rang $0<r'<r$ on a 
$\mu(L)<\mu(E)$ 
(resp. $\frac{\chi(L(n))}{r'}<\frac{\chi(E(n))}{r}$ pour $n\gg 0$). 
On a les implications suivantes :
$$\mu\text{-stable}\Longrightarrow\text{stable}\Longrightarrow
\text{semi-stable}\Longrightarrow\mu\text{-semi-stable}$$
Supposons enfin $\text{Pic}(X)\simeq\mathbb{Z}$ et soit $F$ un faisceau 
r\'eflexif de rang 2 sur $X$, de premi\`ere classe de Chern $c_{1}(F)=0$ ou 
$c_{1}(F)=-1$. Alors $F$ est stable si et seulement si $h^{0}(F)=0$ 
et si $c_{1}(F)=0$ alors $F$ est semi-stable 
si et seulement si $h^{0}(F(-1))=0$ ([H2] lemme 3.1).\\
\newline
\indent (1.3) Soient $X\subset\mathbb{P}^{4}$ une hypersurface cubique 
lisse et $\mathcal{O}_{X}(1)$ le g\'en\'erateur tr\`es ample 
de $\text{Pic}(X)$. Les $\mathbb{Z}$-modules  
$H^{2}(X,\mathbb{Z})$, 
$H^{4}(X,\mathbb{Z})$ et $H^{6}(X,\mathbb{Z})$ 
sont libres de rang 1. On identifie 
ainsi les classes de Chern d'un faisceau coh\'erent sur $X$ \`a des 
entiers relatifs. Nous \'etudions ici l'espace des modules des faisceaux 
semi-stables de rang 2 sur $X$. Nous d\'emontrons le :\\
\newline
\textsc{Th\'eor\`eme 1.4}.$-$\textit{Soient $X\subset\mathbb{P}^{4}$ 
une hypersurface cubique lisse et $B$ la 
surface de Fano de $X$. Alors l'espace des modules $M_{X}$ des faisceaux 
semi-stables de rang 2 sur $X$ de classe de Chern $c_{1}=0$, $c_{2}=2$ et 
$c_{3}=0$ est isomorphe \`a l'\'eclatement d'un translat\'e de la 
surface $-B$ dans la jacobienne interm\'ediaire $J(X)$.}\\
\newline 
\indent (1.5) Soit 
$X$ une vari\'et\'e projective lisse de dimension au moins 2 
et $E$ un fibr\'e vectoriel  
de rang 2 sur $X$. S'il existe une section globale dont le lieu 
des z\'eros $Y$ est de codimension pure 2 alors on a 
une suite exacte ([H1]) :
$$0\longrightarrow\mathcal{O}_{X}\longrightarrow E\longrightarrow
I_{Y}\otimes\text{det}(E)\longrightarrow 0$$
\indent (1.6) \textit{Fibr\'es de rang 2 et construction de Serre}.$-$
Supposons $X$ de dimension au moins 3. Soit $L$ un fibr\'e inversible 
sur $X$ tel que $h^{1}(L^{-1})=0$ et $h^{2}(L^{-2})=0$ et soit 
$Y\subset X$ un sous-sch\'ema ferm\'e de codimension pure 2. On a un 
isomorphisme $\text{Ext}^{1}_{X}(I_{Y}\otimes L,\mathcal{O}_{X})=
H^{0}(\mathcal{O}_{Y}).$ Le sous-sch\'ema $Y$ est 
le lieu des z\'eros d'une section d'un fibr\'e $E$ de rang 2 sur $X$ 
de d\'eterminant $L$ si et seulement si $Y$ est localement 
intersection compl\`ete et $\omega_{Y}=(\omega_{X}\otimes L)_{|Y}$.\\
\newline
\indent (1.7) Soit $X\subset\mathbb{P}^{4}$ une hypersurface cubique lisse. 
Nous montrons que les fibr\'es vectoriels stables sont associ\'es aux 
quintiques elliptiques normales trac\'ees sur $X$ par la construction 
de Serre (2.4), au moyen du :\\   
\newline
\indent (1.8) \textit{Crit\`ere de Mumford-Castelnuovo}.$-$Soit 
$F$ un faisceau coh\'erent sur une vari\'et\'e projective  
$X$ tel que $h^{i}(F(-i))=0$ pour $i\ge 1$. 
Alors $h^{i}(F(k))=0$ pour $i\ge 1$ et $k\ge -i$ et 
$F$ est engendr\'e par ses sections globales ([Mum] lect. 14).\\
\newline
\indent (1.9) Soit $X\subset\mathbb{P}^{4}$ une hypersurface cubique lisse.
Nous montrons que les faisceaux stables non localement 
libres sont param\'etr\'es par les coniques lisses trac\'ees sur $X$ et que 
les faisceaux strictement semi-stables sont param\'etr\'es par les 
couples de droites de $X$ (3.5). Nous montrons enfin que la seconde 
classe de Chern d\'efinit un morphisme vers la jacobienne 
interm\'ediaire $J(X)$. Ce morphisme 
est birationnel ([I-M]) et identifie $M_{X}$ \`a l'\'eclatement d'une 
surface lisse dans $J(X)$ (4.8).\\
\newline
\indent (1.10) Soient $F$ un faisceau coh\'erent sur un sch\'ema $X$ et 
$Y\subset X$ un sous-sch\'ema ferm\'e. La restriction de $F$ \`a $Y$ 
sera not\'ee $F_{Y}$.\\
\newline
\textbf{Remerciements}.$-$Je tiens \`a exprimer toute ma gratitude \`a
Arnaud Beauville pour m'avoir soumis ce probl\`eme et pour l'aide qu'il m'a
apport\'ee. Je remercie \'egalement le $referee$ pour ces nombreuses 
remarques et pour avoir relev\'e une erreur dans la preuve de la proposition 
3.1.\\
\newline
\textbf{2. Fibr\'es de rang 2 stables sur la cubique de 
$\mathbb{P}^{4}$}\\
$\ $
\newline
\textsc{Lemme 2.1}.$-$\textit{Soient $X\subset\mathbb{P}^{N}$ 
une vari\'et\'e de dimension $n\ge 2$ et $E$ un fibr\'e de rang 2 
$\mu$-semi-stable de premi\`ere classe de 
Chern $c_{1}(E)=0$. Si $h^{0}(E)\neq 0$ alors 
le lieu des z\'eros d'une section globale non nulle est de codimension 
pure 2 ou bien ladite section ne s'annule pas et $c_{2}(E)=0$.}\\
\newline
\textit{D\'emonstration}.$-$Le fibr\'e $E$ est de rang 2 et toute 
section non triviale s'annule donc en codimension au plus 2 ou bien 
ne s'annule pas. S'il existe une section partout non nulle alors 
$E$ est extension de $\mathcal{O}_{X}$ par $L$ avec $c_{1}(L)=0$ et 
$c_{2}(E)=0$. Supposons qu'une section de $E$ s'annule en codimension 1
et soit $D$ la partie de codimension 1 du lieu des z\'eros 
de ladite section. On a ainsi $h^{0}(E(-D))\neq 0$ 
et $\mu(\mathcal{O}_{X}(D))\le 0$ puisque $E$ est semi-stable. 
Or $D$ est effectif et on a donc $D=0$.\qed\\ 
\newline
\textsc{Lemme 2.2}.$-$\textit{Soient $S\subset\mathbb{P}^{3}$ une surface 
cubique lisse et $E$ un fibr\'e de rang 2 $\mu$-semi-stable de classes 
de Chern $c_{1}(E)=0$ et $c_{2}(E)=2$. Si $h^{0}(E)=0$ 
alors $h^{1}(E(n))=0$ pour $n\in\mathbb{Z}$ et 
$h^{2}(E(n))=0$ pour $n\ge -1$. Si 
$h^{0}(E)\neq 0$ alors $h^{0}(E)=1$, 
$h^{1}(E(n))=0$ pour $n\le -2$ et $n\ge 1$, 
$h^{1}(E(-1))=h^{1}(E)=1$ et 
$h^{2}(E(n))=0$ pour $n\ge 0$.}\\
\newline
\textit{D\'emonstration}.$-$\textit{Supposons $h^{0}(E)=0$.}$-$On 
a $h^{1}(E)=h^{0}(E)=0$ puisque $h^{2}(E)=h^{0}(E(-1))=0$ et 
$\chi(E)=0$. On a enfin 
$h^{2}(E(-1))=h^{0}(E)=0$. Finalement $h^{i}(E(1-i))=0$ pour $i\ge 1$ 
et le lemme est une cons\'equence de (1.8).\\
\textit{Supposons $h^{0}(E)\neq 0$.}$-$Le fibr\'e $E$ 
est semi-stable et le lieu des z\'eros d'une section globale non nulle 
est donc de codimension pure 2 (2.1). On a donc une suite exacte (1.5) : 
$$
0\longrightarrow\mathcal{O}_{S}\longrightarrow E 
\longrightarrow I_{Z}\longrightarrow 0
$$
o\`u $Z$ est un sous-sch\'ema ferm\'e de dimension 0 et de longueur 2. 
On en d\'eduit en particulier $h^{0}(E)=1$ 
et $h^{1}(E)=1$ puisque $h^{2}(E)=h^{0}(E(-1))=0$ et 
$\chi(E)=0$. 
L'application naturelle  
$H^{0}(\mathcal{O}_{S}(1))\longrightarrow H^{0}(\mathcal{O}_{Z}(1))$ est 
surjective puisque $\ell(Z)=2$ et on a donc $h^{1}(I_{Z}(1))=0$. 
On en d\'eduit $h^{1}(E(1))=0$. Finalement $h^{i}(E(2-i))=0$ pour 
$i\ge 1$ et le lemme est une cons\'equence de (1.8).\qed\\
\newline
\textsc{Lemme 2.3}.$-$\textit{Soient $S\subset\mathbb{P}^{3}$ une surface 
cubique lisse et $E$ un fibr\'e de rang 2 $\mu$-semi-stable de classes 
de Chern $c_{1}(E)=0$ et $c_{2}(E)=1$. Si 
$h^{0}(E)\neq 0$ alors $h^{0}(E)=1$, 
$h^{1}(E(n))=0$ pour $n\in\mathbb{Z}$
et $h^{2}(E(n))=0$ pour $n\ge 0$.}\\
\newline
\textit{D\'emonstration}.$-$Le fibr\'e $E$ 
est semi-stable et le lieu des z\'eros d'une section globale non nulle 
est donc de codimension pure 2 (2.1). On a une suite exacte (1.5) :
$$
0\longrightarrow\mathcal{O}_{S}\longrightarrow E
\longrightarrow I_{Z}\longrightarrow 0
$$
o\`u $Z$ est un point de $S$. On en d\'eduit 
$h^{0}(E)=1$. La suite exacte :
$$
0\longrightarrow\mathcal{O}_{S}(n)\longrightarrow E(n)
\longrightarrow I_{Z}(n)\longrightarrow 0
$$
donne $h^{1}(E(n))=0$ pour 
$n\ge 0$ puisque $h^{1}(\mathcal{O}_{S}(n))=0$ et 
$h^{1}(I_{Z}(n))=0$ pour $n\ge 0$. On en d\'eduit 
$h^{1}(E(n))=h^{1}(E(-n-1))=0$ pour 
$n<0$. On a enfin 
$h^{2}(E(n))=h^{0}(E(-1-n))=0$ pour $n\ge 0$.\qed\\
\newline
\textsc{Th\'eor\`eme 2.4}.$-$\textit{Soient $X\subset\mathbb{P}^{4}$
une cubique lisse et $E$ un fibr\'e de rang 2
stable de classes de Chern $c_{1}(E)=0$ et $c_{2}(E)=2$. Alors $E(1)$ 
est engendr\'e par ses sections globales.}\\
\newline
\textit{D\'emonstration}.$-$Soit $S\in |\mathcal{O}_{X}(1)|$ une section 
hyperplane g\'en\'erique de $X$ tel que le fibr\'e $E_{S}$ soit 
$\mu$-semi-stable relativement \`a 
la polarisation $\mathcal{O}_{S}(1)$ ([M] thm. 3.1).\\
\newline
\textit{Supposons $h^{0}(E_{S})=0$}.$-$Il suffit de 
prouver $h^{i}(E(1-i))=0$ pour $i\ge 1$ (1.8).
Consid\'erons la suite 
exacte :
$$
0\longrightarrow E(n-1) \longrightarrow E(n) \longrightarrow E_{S}(n)
\longrightarrow 0 
$$
On a $h^{1}(E(n))\le h^{1}(E(n-1))$ puisque 
$h^{1}(E_{S}(n))=0$ pour $n\in\mathbb{Z}$ (2.2).
On en d\'eduit $h^{1}(E(n))=0$ pour $n\in\mathbb{Z}$ puisque 
$h^{1}(E(n))=0$ pour $n\ll 0$ puis $h^{2}(E(n))=0$ pour $n\in\mathbb{Z}$. 
On a enfin $h^{3}(E(-2))=h^{0}(E)=0$.\\
\newline
\textit{Supposons $h^{0}(E_{S})\neq 0$ et montrons que nous 
aboutissons \`a une contradiction.}$-$Le 
fibr\'e $E(2)$ 
est alors engendr\'e par ses sections globales. 
Il suffit en effet de prouver
$h^{i}(E(2-i))=0$ pour $i\ge 1$ (1.8). 
Consid\'erons \`a nouveau la suite exacte :
$$
0\longrightarrow E(n-1) \longrightarrow E(n) \longrightarrow E_{S}(n)
\longrightarrow 0
$$
On a $h^{1}(E(n))\le h^{1}(E(n-1))$ pour $n\le-2$ puisque 
$h^{1}(E_{S}(n))=0$ pour $n\le-2$ (2.2).
On en d\'eduit $h^{1}(E(-n))=0$ pour $n\ge 2$ puisque 
$h^{1}(E(n))=0$ pour $n\ll 0$.
Calculons $h^{1}(E(1))$. On a $h^{2}(E)=h^{1}(E(-2))=0$ et 
$h^{3}(E)=h^{0}(E(-2))=0$. Puisque 
$\chi(E)=0$ on a donc $h^{1}(E)=0$ et la suite exacte :
$$
0\longrightarrow E \longrightarrow E(1) \longrightarrow E_{S}(1)
\longrightarrow 0
$$
entra\^{\i}ne $h^{1}(E(1))=h^{1}(E_{S}(1))=0$ (2.2). 
On a enfin $h^{3}(E(-1))=h^{0}(E(-1))=0$. 
Le fibr\'e $E(2)$ est donc engendr\'e par ses sections globales.\\
\indent Si l'une des sections du fibr\'e $E(2)$ est partout non nulle 
alors $E(2)$ est isomorphe au fibr\'e 
$\mathcal{O}_{X}(2)\oplus\mathcal{O}_{X}(-2)$ et $c_{2}(E)=-12$ ce qui 
est absurde. On a donc une suite exacte (1.5) :
$$
0\longrightarrow\mathcal{O}_{X}(-4) \longrightarrow E(-2) 
\longrightarrow I_{C} \longrightarrow 0
$$
o\`u $C\subset X$ est une courbe lisse de degr\'e $c_{2}(E(2))=14$. 
On a $h^{1}(I_{C})=0$ et 
la courbe $C$ est donc connexe. On a  
$\omega_{C}=\mathcal{O}_{C}(2)$ (1.6) et $g(C)=15$. 
Enfin, la courbe $C$ est non d\'eg\'en\'er\'ee puisque le fibr\'e $E$ est 
stable. Calculons $h^{0}(\mathcal{O}_{C}(1))$. La suite exacte :
$$
0\longrightarrow\mathcal{O}_{X}(-3) \longrightarrow E(-1) 
\longrightarrow I_{C}(1) \longrightarrow 0
$$
entra\^{\i}ne l'\'egalit\'e 
$h^{1}(I_{C}(1))=h^{1}(E(-1))$ puisque 
$h^{1}(\mathcal{O}_{X}(-3))=0$ et $h^{2}(\mathcal{O}_{X}(-3))=0$. 
La suite exacte :
$$
0\longrightarrow E(-2) \longrightarrow E(-1) 
\longrightarrow E_{S}(-1) \longrightarrow 0
$$
donne $h^{1}(E(-1))=h^{1}(E_{S}(-1))=1$ (2.2) puisque 
$h^{1}(E(-2))=0$ et $h^{2}(E(-2))=h^{1}(E)=0$. On a 
donc $h^{1}(I_{C}(1))=1$.
On d\'eduit de la suite exacte :
$$
0\longrightarrow I_{C}(1)\longrightarrow\mathcal{O}_{X}(1)
\longrightarrow\mathcal{O}_{C}(1)\longrightarrow 0
$$
que $h^{0}(\mathcal{O}_{C}(1))=6$ puisque 
$h^{0}(I_{C}(1))=0$ et $h^{1}(\mathcal{O}_{X}(1))=0$. 
La courbe $C$ est donc la projection dans $\mathbb{P}^{4}$ 
d'une courbe de Castelnuovo de $\mathbb{P}^{5}$
et le lemme 2.5 fournit la contradiction cherch\'ee.\qed\\
\newline
\textsc{Lemme 2.5}.$-$\textit{Soit $C\subset\mathbb{P}^{5}$ une courbe 
non d\'eg\'en\'er\'ee de genre 15 et de degr\'e 14 (courbe de 
Castelnuovo). Soit $O\in\mathbb{P}^{5}$ ($O\not\in C$) tel 
que la projection \`a partir
de $O$ induise un plongement de $C$ dans $\mathbb{P}^{4}$. L'image 
de $C$ dans $\mathbb{P}^{4}$ n'est alors contenue dans aucune cubique 
lisse.}\\
\newline
\textit{D\'emonstration}.$-$La courbe $C$ est 
contenue dans une surface
irr\'eductible $S\subset\mathbb{P}^{5}$ de degr\'e 4. Ladite surface
$S$ et la courbe $C$ sont ([A-C-G-H]) :\\
ou bien\\
\indent $\bullet$ la surface de Veronese et $C$ est l'image d'une courbe
plane de degr\'e 7 par le plongement de Veronese,\\
ou bien\\
\indent $\bullet$ l'image de 
$S_{2k}=\mathbb{P}_{\mathbb{P}^{1}}
(\mathcal{O}_{\mathbb{P}^{1}}\oplus\mathcal{O}_{\mathbb{P}^{1}}(-2k))$, 
$k\in\{0,1,2\}$, par 
le morphisme $\varphi_{k}$ associ\'e au syst\`eme lin\'eaire 
$|C_{0}+(k+2)f|$ et $C\in |4C_{0}+(4k+6)f|$, o\`u $C_{0}$ est la section 
associ\'ee au fibr\'e naturel $\mathcal{O}_{S_{2k}}(1)$ 
$(C_{0}^{2}=-2k)$ et $f$ une g\'en\'eratrice de la surface r\'egl\'ee 
$S_{2k}$. Pour $k\in\{0,1\}$ le morphisme $\varphi_{k}$ est un
plongement ferm\'e,  
$S=\varphi_{2}(S_{4})$ est un c\^one au-dessus d'une courbe 
rationnelle lisse de degr\'e 4 et le morphisme $\varphi_{2}$ s'identifie
\`a l'\'eclatement de $\varphi_{2}(S_{4})$ en son sommet.\\
\indent Notons $\pi$ la projection consid\'er\'ee
et $\pi(S)$ l'image de $S$ par l'application rationnelle $\pi$.
Si $\pi(S)$ est de dimension 1 alors $S$ est un c\^one au dessus de $C$ 
isomorphe \`a $\varphi_{2}(S_{4})$ ce qui absurde puisque 
$g(C)\ge 1$. La vari\'et\'e $\pi(S)$ est donc de dimension 2. Si $S$ 
est un c\^one alors son sommet et le point de projection 
sont donc distincts.\\
\indent Supposons la courbe $C\subset\mathbb{P}^{4}$ contenue dans une 
cubique lisse $X$ et notons 
$\overline{X}\subset\mathbb{P}^{5}$  le c\^one de sommet $O$ et de 
base $X$.\\
\indent Supposons que la cubique $\overline{X}$ ne contienne pas la 
surface $S$. L'hypersurface $\overline{X}$ d\'ecoupe alors sur $S$ une 
courbe de degr\'e 12 et ne peut donc pas contenir la courbe $C$. 
La cubique $\overline{X}$ contient donc la surface S. On en d\'eduit en
particulier que $\pi(S)\subset X$.\\
\indent Supposons $O\in S$. Si $S$ est l'une des deux surfaces 
$\varphi_{k}(S_{2k})$ avec $k\in\{0,1,2\}$ alors la g\'en\'eratrice $f$ 
passant par $O$ est contract\'ee par $\pi$. Or $C.f=4$ et $\pi$ 
ne peut donc pas induire 
un plongement de $C$ dans $\mathbb{P}^{4}$. 
Si $S$ est la surface de Veronese alors l'application 
rationnelle 
$\mathbb{P}^{2}\dashrightarrow \mathbb{P}^{4}$
obtenue est d\'efinie par le syst\`eme lin\'eaire des coniques
passant par un point. Ce syst\`eme lin\'eaire induit un plongement 
de la surface de Hirzebruch $\mathbb{F}_{1}$ 
dans $\mathbb{P}^{4}$ dont l'image est une surface de degr\'e 3. 
Or $\mathbb{F}_{1}\subset X$ et ladite surface est 
un diviseur de Cartier associ\'e au fibr\'e 
$\mathcal{O}_{X}(l)$ o\`u $l$ est un entier convenable. Son degr\'e est 
donc $3l$. On en d\'eduit que la surface $\pi(S)$ est une 
section hyperplane de $X$, ce qui est absurde.\\ 
\indent Il nous reste \`a traiter le cas o\`u $O\notin S$.
Notons $d$ le degr\'e de $\pi$. La surface $\pi(S)$ est donc 
de degr\'e 
$\frac{4}{d}$. C'est un diviseur de Cartier associ\'e
au fibr\'e $\mathcal{O}_{X}(l)$
o\`u $l$ est un entier convenable. Son degr\'e est donc $3l$ ce qui 
constitue la contradiction cherch\'ee puisque l'\'egalit\'e 
$3ld=4$ est impossible avec $l$ et $d$ entiers.\qed\\
\newline
\textsc{Corollaire 2.6}.$-$\textit{Le fibr\'e $E$ est associ\'e \`a 
une quintique elliptique lisse non d\'eg\'en\'er\'ee par la construction 
de Serre.}\\
\newline
\textit{D\'emonstration}.$-$Il est donn\'e par l'extension (1.5) :
$$
0\longrightarrow\mathcal{O}_{X}(-2) \longrightarrow E(-1)
\longrightarrow I_{C}\longrightarrow 0
$$
o\`u $C$ est une courbe lisse. On a en particulier 
$h^{1}(I_{C})=0$ et la courbe $C$ est donc connexe. 
On a $\omega_{C}=\mathcal{O}_{C}$ (1.6) et la courbe $C$ est donc 
une courbe elliptique de degr\'e $c_{2}(E(1))=5$. Enfin la courbe 
$C$ est lin\'eairement normale puisque $E$ est stable.\qed\\
\newline
\textbf{3. Faisceaux de rang 2 semi-stables sur la cubique 
de $\mathbb{P}^{4}$}\\
\newline
\textsc{Proposition 3.1}.$-$\textit{Soient X une cubique lisse de 
$\mathbb{P}^{4}$ et $E$ un faisceau de rang 2 semi-stable  
de classes de Chern $c_{1}(E)=0$, $c_{2}(E)=2$ et $c_{3}(E)=0$. 
Soit $F$ le bidual de $E$. Alors ou bien $E$ est localement libre ou bien 
$F$ est localement libre de seconde classe de Chern $c_{2}(F)=1$
et $h^{0}(F)=1$ 
ou bien $F=H^{0}(F)\otimes\mathcal{O}_{X}$.}\\
\newline
\textit{D\'emonstration}.$-$Soit  
$S\in|\mathcal{O}_{X}(1)|$ une section hyperplane g\'en\'erique 
telle que $E_{S}$ soit $\mu$-semi-stable relativement \`a la polarisation 
$\mathcal{O}_{S}(1)$ ([M] thm 3.1) et telle que 
$F_{S}$ soit isomorphe au bidual de $E_{S}$. Le faisceau $F$ est $\mu$-semi-stable.
Le faisceau $F_{S}$ est localement libre de rang $2$ et  
$\mu$-semi-stable de premi\`ere classe de Chern $c_{1}(F_{S})=0$ 
([H2]). Notons $R$ le conoyau de l'inclusion canonique $E\subset F$. 
Le faisceau $E$ est sans torsion et $R$ est de dimension au plus 1. 
On a les formules $c_{2}(F_{S})=c_{2}(E_{S})+c_{2}(R_{S})=2-\ell(R_{S})$ 
et $\chi(F_{S})=\ell(R_{S})$. On en d\'eduit la relation  
$h^{0}(F_{S})=h^{1}(F_{S})+\ell(R_{S})$ puisque 
$h^{2}(F_{S})=h^{0}(F_{S}(-1))=0$.
Supposons $h^{0}(F_{S})\ge 1$. Le lieu des z\'eros d'une 
section non nulle est ou bien vide ou bien de codimension pure 2 (2.1). 
S'il est vide alors le fibr\'e $F_{S}$ est trivial et s'il est de 
codimension pure 2 alors $h^{0}(F_{S})=1$. On a donc  
$\ell(R_{S})\in\{0,1,2\}$ et $c_{2}(F_{S})\in\{0,1,2\}$.\\
\indent Consid\'erons la suite exacte de restriction \`a une section 
hyperplane :
$$
0\longrightarrow F(n-1)\longrightarrow F(n)\longrightarrow
F_{S}(n)\longrightarrow 0
$$
On a $h^{1}(F_{S}(n))=0$ pour $n\le -2$ et $n\ge 1$ (2.2 et 2.3). On en 
d\'eduit $h^{1}(F(n-1))\ge h^{1}(F(n))$ pour $n\le-2$ et 
$h^{2}(F(n-1))\le h^{2}(F(n))$ pour $n\ge 1$. Or 
$h^{1}(F(n))=0$ pour $n\ll 0$ ([H2] thm. 2.5) et 
$h^{2}(F(n))=0$ pour $n\gg 0$ et on a donc 
$h^{1}(F(n))=0$ pour $n\le -2$
$h^{2}(F(n))=0$ pour $n\ge 0$. 
On a enfin $h^{3}(F)=h^{0}(F^{*}(-2))=h^{0}(F(-2))=0$ ([H2] prop. 1.10) 
et $h^{0}(F)\le h^{0}(F_{S})$.\\
\newline
\textit{Supposons $\ell(R_{S})=0$}$-$Alors $c_{2}(F)=2$  
et $\chi(F)=\frac{c_{3}(F)}{2}$.
On en d\'eduit la formule $\frac{c_{3}(F)}{2}=h^{0}(F)-h^{1}(F)$. Or 
$c_{3}(F)\ge 0$ ([H2] prop. 2.6) et $h^{0}(F)\le h^{0}(F_{S})\le 1$
et on a donc $c_{3}(F)=0$ ou $2$.\\
\indent Si $c_{3}(F)=0$ alors les faisceaux  
$E$ et $F$ sont canoniquement isomorphes et localement libres 
([H2] prop. 2.6).\\ 
\indent Si $c_{3}(F)=2$ alors $h^{0}(F)=1$ et $h^{1}(F)=0$. Le 
faisceau $R$ est donc de dimension 0 et 
$\ell(R)=\chi(F)-\chi(E)=1$. On a donc $R=k(p)$ avec $p\in X$. Puisque 
$h^{0}(F)=1$ on a un morphisme non nul 
$\mathcal{O}_{X}\longrightarrow F$.
De plus,
$\chi(E(n))=n^{3}+3n^{2}+2n$ et 
$\chi(\mathcal{O}_{X}(n))=\frac{n^{3}}{2}+\frac{3n^{2}}{2}+2n+1$ et 
on en d\'eduit $h^{0}(E)=0$ puisque $E$ est semi-stable.  
Le morphisme induit $\mathcal{O}_{X}\longrightarrow R$ est donc 
non nul. Il est surjectif et induit 
une inclusion $I_{p}\subset X$. Or 
$\chi(I_{p}(n))=\frac{n^{3}}{2}+\frac{3n^{2}}{2}+2n$ ce qui est en 
contradiction avec la semi-stabilit\'e de $E$.\\
\newline
\textit{Supposons $\ell(R_{S})\ge 1$.} 
On a donc $h^{0}(F_{S})\ge 1$. Le lieu des z\'eros d'une section globale
non nulle est ou bien vide ou bien de codimension pure 2 (2.1).\\
\newline
\indent Supposons qu'il existe une section non nulle de $F_{S}$ 
dont le lieu des z\'eros est de codimension pure 2. 
Alors $h^{0}(F_{S})=1$. On en d\'eduit 
$h^{1}(F_{S})=0$ puis $\ell(R_{S})=1$ et $c_{2}(F)=1$. 
On en d\'eduit l'in\'egalit\'e 
$\chi(F)=h^{0}(F)-h^{1}(F)=1+\frac{c_{3}(F)}{2}
\le 1-h^{1}(F)$. Puis $c_{3}(F)=0$ puisque $c_{3}(F)\ge 0$ ([H2] prop. 2.6).
Le faisceau $F$ est donc localement libre ([H2] prop. 2.6) de seconde 
classe de Chern $c_{2}(F)=1$ et $h^{0}(F)=1$.\\
\newline
\indent Supposons enfin qu'il existe une section du fibr\'e $F_{S}$ 
ne s'annulant pas auquel cas ledit fibr\'e est isomorphe au fibr\'e 
$H^{0}(F_{S})\otimes\mathcal{O}_{S}$ et donc 
$\ell(R_{S})=2$ et $c_{2}(F)=0$.
On en d\'eduit l'in\'egalit\'e
$\chi(F)=\frac{c_{3}(F)}{2}+2=h^{0}(F)-h^{1}(F)\le 2-h^{1}(F)$. 
puis $c_{3}(F)=0$ puisque $c_{3}(F)\ge 0$ ([H2] prop. 2.6)
et $h^{0}(F)=2$. Le faisceau $F$ est donc localement libre ([H2] prop. 2.6).
Supposons qu'il existe une section globale non nulle de $F$ dont le lieu des 
z\'eros $Z$ est non vide. Le sch\'ema $Z$ est de 
dimension pure 1 puisque 
$h^{0}(F(-1))=0$ et $F$ est donc extension de 
$I_{Z}$ par $\mathcal{O}_{X}$. On en d\'eduit $h^{0}(F)=1$ 
ce qui est absurde. Le faisceau $F$ est donc 
isomorphe au fibr\'e $H^{0}(F)\otimes\mathcal{O}_{X}$.\qed\\
\newline
\textsc{Lemme 3.2}.$-$\textit{Soit $R$ un faisceau 
coh\'erent sur $\mathbb{P}^{n}$ $(n\ge1)$ tel que 
$h^{0}(R(-1))=0$ et 
$\chi(R(n))=n+1$. Il existe alors une 
droite $\ell\subset\mathbb{P}^{n}$ telle que $R=\mathcal{O}_{\ell}$.}\\
\newline
\textit{D\'emonstration}.$-$Le faisceau $R$ est de dimension 1 et 
on a donc $h^{0}(R)=h^{1}(R)+1\ge 1$. Soient $s\in H^{0}(R)$ une 
section non nulle et $I_{Z}$ le noyau de l'application induite 
$\mathcal{O}_{\mathbb{P}^{n}}\longrightarrow R$.
On a $h^{0}(\mathcal{O}_{Z}(-1))=0$ et $Z$ 
est donc de dimension pure 1.
Consid\'erons une section hyperplane g\'en\'erique 
$S\in|\mathcal{O}_{\mathbb{P}^{n}}(1)|$. On a $\ell(R_{S})=1$ et 
l'inclusion $\mathcal{O}_{Z\cap S}\subset R_{S}$ est donc un 
isomorphisme. On en d\'eduit que  
$Z_{\text{red}}$ est une droite $\ell\subset \mathbb{P}^{n}$
et 
que le sch\'ema $Z$ est g\'en\'eriquement r\'eduit le long 
de $\ell$. Le noyau 
de l'application surjective $\mathcal{O}_{Z}
\longrightarrow \mathcal{O}_{\ell}$ est de dimension z\'ero et 
donc trivial puisque $h^{0}(\mathcal{O}_{Z}(-1))=0$. 
On a donc $R=\mathcal{O}_{\ell}$ puisque ces deux faisceaux ont m\^eme 
polyn\^ome caract\'eristique.\qed\\
\newline
\textsc{Lemme 3.3}.$-$\textit{Soit $R$ un faisceau coh\'erent 
sur $\mathbb{P}^{n}$ $(n\ge 1)$ tel que 
$h^{0}(R(-1))=0$
et $\chi(R(n))=2n+2$. Alors il existe deux droites 
$\ell_{1}\subset\mathbb{P}^{n}$ et 
$\ell_{2}\subset\mathbb{P}^{n}$ telles que 
$R$ soit extension de $\mathcal{O}_{\ell_{2}}$ par 
$\mathcal{O}_{\ell_{1}}$ ou bien $R(-1)$ est une 
th\^eta-caract\'eristique sur une conique lisse $C\subset\mathbb{P}^{n}$.}\\
\newline
\textit{D\'emonstration}.$-$Le faisceau $R$ est de dimension 1 et 
on a donc  
$h^{0}(R)=h^{1}(R)+2\ge 2$. 
Soient $s\in H^{0}(R)$ une section non nulle 
et $I_{Z}$ le noyau de l'application induite 
$\mathcal{O}_{\mathbb{P}^{n}}\longrightarrow R$. On a  
$h^{0}(\mathcal{O}_{Z}(-1))=0$ et $Z$ est donc de 
dimension pure 1. Soit $S\in|\mathcal{O}_{\mathbb{P}^{n}}(1)|$ 
une section hyperplane g\'en\'erique.
On a $0<\ell(Z\cap S)\le \ell(R_{S})=2$. Notons $Q$ le conoyau 
de l'inclusion $\mathcal{O}_{Z}\subset R$.\\
\newline
\textit{Supposons $\ell(Z\cap S)=1$}.$-$Le support du sch\'ema $Z$ 
est alors une droite $\ell_{1}$ et ledit sch\'ema est 
g\'en\'eriquement r\'eduit le long de 
$\ell_{1}$. On a donc une application surjective 
$\mathcal{O}_{Z}\longrightarrow\mathcal{O}_{\ell_{1}}$ dont le noyau est  
de dimension z\'ero. Ledit noyau est en fait trivial puisque 
$h^{0}(\mathcal{O}_{Z}(-1))=0$. Enfin, on a 
$\chi(Q(n))=n+1$ et $h^{0}(Q(-1))=0$ 
et le lemme 3.2 permet de conclure.\\
\newline
\textit{Supposons $\ell(Z\cap S)=2$ et $Q$ non trivial}.$-$On a 
$h^{1}(R(-1))=0$ et on a donc $h^{1}(R(k))=0$ pour 
$k\ge -1$ (1.8). On en d\'eduit en particulier $h^{0}(R)=2$ 
et $h^{0}(R(1))=4$. Consid\'erons la suite exacte :
$$
0\longrightarrow\mathcal{O}_{Z}(1)\longrightarrow R(1)\longrightarrow
Q(1)\longrightarrow 0
$$
o\`u $Q$ est de dimension z\'ero. 
Le faisceau $R(1)$ est engendr\'e par ses sections globales (1.8) 
et l'application 
$H^{0}(R(1))\longrightarrow H^{0}(Q(1))$ 
n'est donc pas identiquement nulle. On en d\'eduit 
$h^{0}(\mathcal{O}_{Z}(1))\le 3$ et 
$h^{0}(I_{Z}(1))\ge n-2$. Il existe donc un plan 
$\mathbb{P}^{2}\subset\mathbb{P}^{n}$ contenant le sch\'ema $Z$. Notons 
$J_{Z}$ l'id\'eal de $Z$ dans ledit plan. On a $c_{1}(J_{Z})=-2$ et on 
a donc une inclusion $J_{Z}\subset\mathcal{O}_{\mathbb{P}^{2}}(-2)$ qui 
induit une application surjective 
$\mathcal{O}_{Z}\longrightarrow\mathcal{O}_{C}$ o\`u $C$ est une conique. 
Son noyau est de dimension z\'ero et donc trivial 
puisque $h^{0}(\mathcal{O}_{Z}(-1))=0$. On a donc une suite exacte :
$$0\longrightarrow\mathcal{O}_{C}\longrightarrow R\longrightarrow 
k(p)\longrightarrow 0$$
Si $p\notin C$ alors l'extension pr\'ec\'edente est triviale ce qui est 
absurde puisque $h^{0}(R(-1))=0$. On a donc $p\in C$.
Montrons que $R$ est un $\mathcal{O}_{C}$-
module. Soit $f\in H^{0}(I_{C}(k))$ $(k\ge 0)$ 
l'\'equation d'une hypersurface de degr\'e $k$  
contenant $C$. Consid\'erons le diagramme commutatif suivant :
\begin{equation*}
\begin{CD}
 & & 0 & & 0 & & 0 \\
 & & @VVV @VVV @VVV \\
0 @))) \mathcal{O}_{C}(-k) @))) K(-k) @))) k(p)(-k) \\
 & & @| @VVV @| \\
0 @))) \mathcal{O}_{C}(-k) @))) R(-k)  @))) k(p)(-k) @))) 0 \\
 & & @VV{\times f}V @VV{\times f}V @VV{\times f}V \\
0 @))) \mathcal{O}_{C} @))) R @))) k(p) @))) 0
\end{CD}
\end{equation*}
o\`u les complexes horizontaux sont exacts. Si l'application 
$K\longrightarrow k(p)$ est nulle alors on a une inclusion 
$R(-k)/K(-k)\subset R$ avec $R(-k)/K(-k)$ de dimension 
z\'ero ce qui est impossible puisque 
$h^{0}(R(-1))=0$. Ladite application est donc surjective 
et on en d\'eduit que l'application $R(-k)\longrightarrow R$ est nulle. 
Le faisceau $R$ est donc un $\mathcal{O}_{C}$-module.
On v\'erifie alors qu'on a une suite exacte :
$$0\longrightarrow I_{p}\longrightarrow H^{0}(R)\otimes\mathcal{O}_{C}
\longrightarrow R\longrightarrow 0$$
o\`u $I_{p}$ est l'id\'eal de $p$ dans $C$. Si $C$ est une conique 
lisse alors $R(-1)$ est la th\^eta-caract\'eristique sur $C$. Supposons 
$C$ non lisse et soit $\ell\subset C$ une droite contenant $p$. 
L'inclusion $I_{\ell}\subset H^{0}(R)\otimes\mathcal{O}_{C}$ se factorise
\`a travers l'inclusion $I_{\ell}\subset\mathcal{O}_{C}$ et on 
obtient ainsi une inclusion $\mathcal{O}_{\ell}\subset R$ dont le 
conoyau est \'egalement isomorphe au faisceau structural d'une droite
(3.2).\\
\newline
\textit{Supposons $\ell(Z\cap S)=2$ et $Q=0$}.$-$Le sch\'ema 
$Z_{\text{red}}$ est de degr\'e au plus 2. S'il est 
de degr\'e 2 ledit sch\'ema est ou bien r\'eunion 
de deux droites distinctes ou 
bien une conique lisse. On a alors une application surjective 
$\mathcal{O}_{Z}\longrightarrow\mathcal{O}_{C}$ dont le noyau est 
trivial si $C$ est r\'eunion de droites disjointes et 
support\'e en un point sinon. Ce dernier cas est impossible puisque 
$h^{0}(\mathcal{O}_{Z}(-1))=0$. Si 
$Z_{\text{red}}$ est de degr\'e 1 alors $Z_{\text{red}}$ est une droite 
$\ell\subset \mathbb{P}^{n}$ et on a une surjection 
$\mathcal{O}_{Z}\longrightarrow\mathcal{O}_{\ell}$ dont 
le noyau $K$ v\'erifie $\chi(K(n))=n+1$ et 
$h^{0}(K(-1))=0$. Ce noyau 
est donc isomorphe au faisceau $\mathcal{O}_{\ell}$ (3.2).\qed\\
\newline
\textsc{Lemme 3.4}.$-$\textit{Soient $X\subset\mathbb{P}^{4}$ une cubique 
lisse et $\theta$ une th\^eta-caract\'eristique sur une conique lisse 
$C\subset X$. On consid\'ere le 
faisceau $E$ noyau de l'application surjective 
$H^{0}(\theta(1))\otimes\mathcal{O}_{X}\longrightarrow\theta(1)$. Alors 
$E$ est stable de classes de Chern $c_{1}(E)=0$, $c_{2}(E)=2$ et 
$c_{3}(E)=0$.}\\
\newline
\textit{D\'emonstration}.$-$Le calcul des classes de Chern de $E$ est 
imm\'ediat. Soit $F\subset E$ un sous-faisceau 
de rang 1 de $E$. Le faisceau $F$ est de la forme 
$I_{Z}(a)$ o\`u $Z\subset X$ est un sous-sch\'ema ferm\'e de 
dimension au plus 1 et 
$a\in\mathbb{Z}$. On a un diagramme commutatif \`a lignes et 
colonnes exactes : 
\begin{equation*}
\begin{CD}
 & & 0 & & 0 & & & & \\
 & & @VVV @VVV & & & & \\
 & & H^{0}(\theta(1))\otimes I_{C} @= 
 H^{0}(\theta(1))\otimes I_{C} & & & & \\
 & & @VVV @VVV & & & & \\
0 @))) E @))) H^{0}(\theta(1))\otimes\mathcal{O}_{X} @))) 
 \theta(1) @))) 0 \\
 & & @VVV @VVV @| & & \\
0 @))) \theta @))) H^{0}(\theta(1))\otimes\mathcal{O}_{C} 
 @))) \theta(1) @))) 0 \\
 & & @VVV @VVV & & & & \\
 & & 0 & & 0 & & & & \\
\end{CD}
\end{equation*}
Notons $F_{0}$ le noyau de l'application induite 
$F\longrightarrow\theta$. On a une 
inclusion $F_{0}\subset H^{0}(\theta(1))\otimes I_{C}$. Le faisceau 
$H^{0}(\theta(1))\otimes I_{C}$ est $\mu$-semi-stable de pente nulle 
et on a donc 
$c_{1}(F)=c_{1}(F_{0})\le c_{1}(H^{0}(\theta(1))\otimes I_{C})=0$ 
puisque $\theta$ est de dimension 1.\\
\indent Si $c_{1}(F)<0$ on a $\chi(F(n))<\frac{1}{2}\chi(E(n))$ 
pour $n\gg 0$ par un calcul classique. Si 
$c_{1}(F)=0$ alors $F=I_{Z}$ avec $\text{codim}(Z)\ge 2$ et on a donc 
$I_{Z}^{**}=\mathcal{O}_{X}$. L'inclusion 
$I_{Z}\subset H^{0}(\theta(1))\otimes\mathcal{O}_{X}$ d\'eduite de 
l'inclusion $E\subset H^{0}(\theta(1))\otimes\mathcal{O}_{X}$ est donc 
donn\'ee par une section non nulle $s\in H^{0}(\theta(1))$. 
L'application induite 
$I_{Z}\longrightarrow\theta(1)$ associe donc la section 
$f_{|C}s$ \`a la fonction $f$. La section $s$ \'etant 
g\'en\'eriquement non nulle on en d\'eduit $I_{Z}\subset I_{C}$ et 
donc $\chi(I_{Z}(n))\le\chi(I_{C}(n))<\frac{1}{2}\chi(E(n))$ 
pour $n\gg 0$ puisque $\chi(I_{C}(n))=\frac{n^{3}}{2}+\frac{3n^{2}}{2}$
et $\frac{\chi(E(n))}{2}=\frac{n^{3}}{2}+\frac{3n^{2}}{2}+n$.\qed\\
\newline
\textsc{Th\'eor\`eme 3.5}.$-$\textit{Soient X une cubique lisse de 
$\mathbb{P}^{4}$ et $E$ un faisceau de rang 2 semi-stable  
de classes de Chern $c_{1}(E)=0$, $c_{2}(E)=2$ et $c_{3}(E)=0$. Si 
$E$ est stable alors ou bien $E$ est localement 
libre ou bien $E$ est associ\'e 
\`a une conique lisse $C\subset X$ (3.4). Si 
$E$ est semi-stable non stable alors 
le gradu\'e de $E$ est somme directe des id\'eaux de deux  
droites de $X$.}\\
\newline
\textit{D\'emonstration}.$-$Soit $F$ le bidual de $E$ et $R$ 
le conoyau de l'injection canonique $E\subset F$. Le 
faisceau $E$ est localement libre 
ou bien $F$ est localement libre de seconde classe de Chern $c_{2}(F)=1$ 
et $h^{0}(F)=1$ ou bien 
$F=H^{0}(F)\otimes\mathcal{O}_{X}$ (3.1). On a 
$\chi(E(n))=n^{3}+3n^{2}+2n$ et 
$\chi(\mathcal{O}_{X}(n))=\frac{n^{3}}{2}+\frac{3n^{2}}{2}+2n+1$ et 
on en d\'eduit $h^{0}(E)=0$ puisque $E$ est semi-stable.\\
\newline
\textit{Supposons $E$ localement libre}.$-$On a 
$h^{0}(E)=0$ puisque $E$ est semi-stable 
et le fibr\'e $E$ est donc stable.\\
\newline
\textit{Supposons $F$ localement libre de seconde classe de Chern 
$c_{2}(F)=1$ et $h^{0}(F)=1$}.$-$
Alors $\chi(R(n))=n+1$. Soit $s\in H^{0}(F)$ une section 
non nulle. Elle s'annule le long 
d'une droite ${\ell}_{2}\subset X$ (2.1). 
On a $h^{0}(E)=0$ et la section $s$ de $F$ induit 
une application non nulle 
$\mathcal{O}_{X}\longrightarrow R$.
Notons $I_{Z}$ le noyau de l'application pr\'ec\'edente. 
Le sch\'ema $Z$ est de dimension 1. Sinon on aurait
une inclusion $I_{Z}\subset E$ avec 
$\chi(I_{Z}(n))=\frac{n^{3}}{2}+\frac{3n^{2}}{2}+2n+1-\ell(Z)$ ce qui 
est impossible par semi-stabilit\'e de $E$. 
Soit $S\in|\mathcal{O}_{X}(1)|$ une section hyperplane g\'en\'erique. 
L'inclusion $\mathcal{O}_{Z\cap S}\subset R_{S}$ est un isomorphisme 
puisque $\ell(R_{S})=1$. Aussi 
la composante de dimension 1 du support de $Z$ est une 
droite $\ell_{1}$ et on a donc une application surjective 
$\mathcal{O}_{Z}\longrightarrow\mathcal{O}_{\ell_{1}}$. 
On en d\'eduit $R=\mathcal{O}_{{\ell}_{1}}$ puisque ces deux faisceaux 
ont m\^eme polyn\^ome caract\'eristique.
L'application $\mathcal{O}_{X}\longrightarrow\mathcal{O}_{{\ell}_{1}}$ 
est non nulle et les droites ${\ell}_{1}$ et ${\ell}_{2}$ sont donc
disjointes. Consid\'erons le diagramme commutatif \`a lignes et 
colonnes exactes :
\begin{equation*}
\begin{CD}
 & & 0 & & 0 \\
 & & @VVV @VVV & & & & & \\
0 @))) I_{{\ell}_{1}} @))) \mathcal{O}_{X} 
@))) \mathcal{O}_{{\ell}_{1}} @))) 0 \\
 & &  @VVV   @VVV   @| \\
0 @))) E @))) F @))) \mathcal{O}_{{\ell}_{1}} @))) 0 \\
 & & @VVV @VVV & & & & & \\
 & & I_{\ell_{2}} @= I_{\ell_{2}} & & & & & \\
 & & @VVV @VVV & & & & & \\
 & &  0 & & 0 & & & & & \\
\end{CD}
\end{equation*}
On en d\'eduit que le faisceau $E$ est semi-stable 
non stable et que son gradu\'e est 
$I_{{\ell}_{1}}\oplus I_{{\ell}_{2}}$.\\
\newline
\textit{Supposons $F=H^{0}(F)\otimes\mathcal{O}_{X}$}.$-$
On a 
$\chi(R(n))=2n+2$
et $h^{0}(E)=0$ puisque $E$ est semi-stable. On en d\'eduit 
en particulier $h^{0}(R)\ge 2$.
Consid\'erons une section hyperplane g\'en\'erique 
$S\in|\mathcal{O}_{X}(1)|$. On a la suite exacte :
$$
0\longrightarrow E_{S} \longrightarrow 
H^{0}(X,F)\otimes\mathcal{O}_{S}
\longrightarrow R_{S} \longrightarrow 0
$$
L'application $H^{0}(H^{0}(X,F)\otimes\mathcal{O}_{S})
\longrightarrow H^{0}(R_{S})$ n'est pas nulle puisque le morphisme 
$H^{0}(F)\otimes\mathcal{O}_{S}\longrightarrow R_{S}$ est surjectif. 
On a donc $h^{0}(E_{S})\le 1$. Si $h^{0}(E_{S})=0$ 
alors l'application $H^{0}(F_{S})
\longrightarrow H^{0}(R_{S})$ est un isomorphisme et l'application 
$H^{0}(F_{S})\otimes H^{0}(\mathcal{O}_{S}(n))
\longrightarrow H^{0}(R_{S}(n))$ est surjective pour tout $n\ge 0$ 
puisqu'il existe une section hyperplane de $S$ \'evitant le support 
de $R_{S}$. Si $h^{0}(E_{S})=1$ alors le quotient 
$H^{0}(F_{S})/H^{0}(E_{S})$ est de dimension 1 et  
on a une application surjective 
$(H^{0}(F_{S})/H^{0}(E_{S}))\otimes
\mathcal{O}_{S}\longrightarrow R_{S}$.
Or $\ell(R_{S})=2$ et 
$\mathcal{O}_{S}(1)$ est tr\`es ample et l'application 
$(H^{0}(F_{S})/H^{0}(E_{S}))\otimes H^{0}(\mathcal{O}_{S}(n))
\longrightarrow H^{0}(R_{S}(n))$
est donc surjective pour $n\ge 1$. Il en r\'esulte que l'application 
$H^{0}(F_{S})\otimes H^{0}(\mathcal{O}_{S}(n))
\longrightarrow H^{0}(R_{S}(n))$ est \'egalement surjective. 
On a donc finalement $h^{1}(E_{S}(n))=0$ pour $n\ge 1$ puisque 
$h^{1}(\mathcal{O}_{S}(n))=0$ pour $n\ge 1$. La suite exacte :
$$
0\longrightarrow E(n-1)\longrightarrow E(n)\longrightarrow E_{S}(n)
\longrightarrow 0
$$
donne $h^{2}(E(n-1)\le h^{2}(E(n))$ pour $n\ge 1$. On a donc 
$h^{2}(E(n))=0$ pour $n\ge 0$ puisque $h^{2}(E(n))=0$ 
pour $n\gg 0$. En particulier $h^{2}(E)=0$. 
On d\'eduit de la suite exacte :
$$
0\longrightarrow E\longrightarrow H^{0}(X,F)\otimes\mathcal{O}_{X}
\longrightarrow R\longrightarrow 0
$$
l'\'egalit\'e $h^{3}(E)=0$. Mais $\chi(E)=0$ et on a donc $h^{1}(E)=0$.
On en d\'eduit $h^{0}(R)=2$ et l'inclusion 
$H^{0}(F)\subset H^{0}(R)$ est donc un isomorphisme.
Montrons alors que l'application de restriction  
$H^{0}(R)\longrightarrow H^{0}(R_{S})$
est injective. Supposons 
le contraire. Il existe donc une section $s\in H^{0}(R)$ non nulle 
dont l'image dans $H^{0}(R_{S})$ est nulle. Notons 
$I_{Z}$ le noyau de l'application 
$\mathcal{O}_{X}\longrightarrow R$ d\'efinie par le section $s$ et 
$Q$ le conoyau de l'inclusion $\mathcal{O}_{Z}\subset R$.
Par hypoth\`ese, l'application 
$\mathcal{O}_{Z\cap S}\longrightarrow R_{S}$ 
est nulle et on a donc $R_{S}=Q_{S}$. Le faisceau $Q$ est donc de 
dimension 1 et $c_{2}(Q)=c_{2}(R)$.
Le sch\'ema $Z$ est donc de 
dimension 0. Or 
$\chi(I_{Z}(n))
=\frac{n^{3}}{2}+\frac{3n^{2}}{2}+2n+1-\ell(Z)$
ce qui est en contradiction avec la semi-stabilit\'e de $E$ puisqu'on 
a une inclusion $I_{Z}\subset E$.
L'application de restriction  
$H^{0}(R)\longrightarrow H^{0}(R_{S})$
est injective et on a donc $h^{0}(R(-1))=0$.
Le faisceau $R(-1)$ est donc ou bien une th\^eta-caract\'eristique 
sur une conique lisse 
$C\subset X$ auquel cas $E$ est stable (3.4) 
ou bien il existe deux droites $\ell_{1}\subset X$ 
et $\ell_{2}\subset X$ telles que $R$ soit extension de 
$\mathcal{O}_{\ell_{1}}$ par $\mathcal{O}_{\ell_{2}}$ (3.3) auquel 
cas on a un diagramme commutatif \`a lignes et colonnes exactes :
\begin{equation*}
\begin{CD}
 & & 0 & & 0 & & 0 \\
 & & @VVV @VVV @VVV \\
 0 @))) I_{\ell_{1}} @))) E @))) I_{\ell_{2}} @))) 0 \\
 & & @VVV @VVV @VVV \\
 0 @))) \mathcal{O}_{X} @))) H^{0}(F)\otimes\mathcal{O}_{X} 
 @))) \mathcal{O}_{X} @))) 0 \\
 & & @VVV @VVV @VVV \\
 0 @))) \mathcal{O}_{\ell_{1}} @))) R @))) 
 \mathcal{O}_{\ell_{2}} @))) 0 \\
 & & @VVV @VVV @VVV \\ 
 & & 0 & & 0 & & 0 \\ 
\end{CD}
\end{equation*}
On en d\'eduit que le gradu\'e de $E$ est le faisceau 
$I_{\ell_{1}}\oplus I_{\ell_{2}}$.\qed\\
\newline
\textbf{4. Espace des modules des faisceaux semi-stables sur la cubique de 
$\mathbb{P}^{4}$}\\
\newline
\indent (4.1) Soit $X\subset\mathbb{P}^{4}$ une hypersurface 
cubique lisse et soit $(Def(E),0)$ l'espace des d\'eformations 
verselles d'un faisceau coh\'erent $E$ sur $X$. 
L'espace tangent \`a $Def(E)$ en $0$ s'identifie \`a l'espace vectoriel 
$\text{Ext}^{1}_{X}(E,E)$. Le germe analytique $Def(E)$ est 
lisse en $0$ si $\text{Ext}^{2}_{X}(E,E)=0$.\\
\newline
\textsc{Lemme 4.2}.$-$\textit{Soient $X\subset\mathbb{P}^{4}$ 
une cubique lisse et $\theta$ une th\^eta-carat\'eristique sur une 
conique lisse $C\subset X$. Soit $E$ le noyau de la surjection 
$H^{0}(\theta(1))\otimes\mathcal{O}_{X}\longrightarrow\theta(1)$. Alors 
$\text{Ext}^{2}_{X}(E,E)$ est nul et l'espace vectoriel complexe 
$\text{Ext}^{1}_{X}(E,E)$ est 
de dimension 5.}\\
\newline
\textit{D\'emonstration}.$-$Soit $F$ le noyau de la surjection 
$\mathcal{O}_{\mathbb{P}^{4}}\oplus\mathcal{O}_{\mathbb{P}^{4}}
\longrightarrow\theta(1)$. On 
v\'erifie que le faisceau $F(1)$ est engendr\'e par ses sections
globales en utilisant le crit\`ere de Mumford-Castelnuovo (1.8). 
On en d\'eduit que le faisceau $E(1)$ est \'egalement engendr\'e par ses 
sections globales puisqu'on a un morphisme surjectif 
$F_{|X}(1)\longrightarrow E(1)$. On a donc 
$\text{Hom}_{X}(E,\theta(-1))\subset 
\text{Hom}_{X}(H^{0}(E(1))\otimes
\mathcal{O}_{X}(-1),\theta(-1))=0$. 
Consid\'erons la suite exacte :
$$\text{Ext}^{2}_{X}(H^{0}(\theta(1))\otimes\mathcal{O}_{X},E)
\longrightarrow \text{Ext}^{2}_{X}(E,E)\longrightarrow 
\text{Ext}_{X}^{3}(\theta(1),E)$$
On a $\text{Ext}^{2}_{X}(H^{0}(\theta(1))\otimes\mathcal{O}_{X},E)
\simeq H^{0}(\theta(1))\otimes H^{2}(E)=0$ et 
$\text{Ext}_{X}^{3}(\theta(1),E)\simeq\text{Hom}_{X}(E,\theta(-1))^{*}=0$ 
et donc $\text{Ext}^{2}_{X}(E,E)=0$. Enfin, 
$\text{Ext}_{X}^{3}(E,E)\simeq\text{Hom}_{X}(E,E(-2))^{*}=0$
et $\text{Hom}_{X}(E,E)\simeq\mathbb{C}$.
L'espace vectoriel complexe $\text{Ext}_{X}^{1}(E,E)$ est donc 
de dimension 5 puisque 
$\chi(E,E)=\sum (-1)^{i}\text{Ext}^{i}_{X}(E,E)=-4$.\qed\\
\newline
\textsc{Lemme 4.3}.$-$\textit{Soit $X\subset\mathbb{P}^{4}$ une cubique 
lisse et soient $\ell_{1}\subset X$ et 
$\ell_{2}\subset X$ deux droites. Le groupe 
$\text{Ext}^{2}_{X}(I_{\ell_{1}},I_{\ell_{2}})$ est nul
et l'espace vectoriel complexe 
$\text{Ext}^{1}_{X}(I_{\ell_{1}},I_{\ell_{2}})$ est 
de dimension 1 si $\ell_{1}\neq\ell_{2}$ et de dimension 2 si
$\ell_{1}=\ell_{2}$.}\\
\newline
\textit{D\'emonstration}.$-$On a un isomorphisme 
$\text{Ext}^{3}_{X}(\mathcal{O}_{\ell_{1}},I_{\ell_{2}})\simeq
\text{Hom}_{X}(I_{\ell_{2}},\mathcal{O}_{\ell_{1}}(-2))^{*}$. Or 
le faisceau $I_{\ell_{2}}(1)$ est engendr\'e par ses sections 
globales et on a donc   
$\text{Hom}_{X}(I_{\ell_{2}},\mathcal{O}_{\ell_{1}}(-2))
\subset\text{Hom}_{X}(H^{0}(I_{\ell_{2}}(1))\otimes\mathcal{O}_{X}(-1),
\mathcal{O}_{\ell_{1}}(-2))=0$. Consid\'erons la suite exacte :
$$\text{Ext}^{2}_{X}(\mathcal{O}_{X},I_{\ell_{2}})
\longrightarrow \text{Ext}^{2}_{X}(I_{\ell_{1}},I_{\ell_{2}})
\longrightarrow \text{Ext}^{3}_{X}(\mathcal{O}_{\ell_{1}},
I_{\ell_{2}})$$
On en d\'eduit $\text{Ext}^{2}_{X}(I_{\ell_{1}},I_{\ell_{2}})=0$
puisque $\text{Ext}^{2}_{X}(\mathcal{O}_{X},I_{\ell_{2}})=0$. 
On a $\text{Ext}^{3}_{X}(I_{\ell_{1}},I_{\ell_{2}})\simeq
\text{Hom}_{X}(I_{\ell_{2}},I_{\ell_{1}}(-2))^{*}=0$ et 
$\chi(I_{\ell_{1}},I_{\ell_{2}})=\sum (-1)^{i}
\text{Ext}^{i}_{X}(I_{\ell_{1}},I_{\ell_{2}})=
\text{Ext}^{0}_{X}(I_{\ell_{1}},I_{\ell_{2}})
-\text{Ext}^{1}_{X}(I_{\ell_{1}},I_{\ell_{2}})
=\chi(I_{\ell_{1}},I_{\ell_{1}})=-1$.
L'espace vectoriel complexe 
$\text{Ext}^{1}_{X}(I_{\ell_{1}},I_{\ell_{2}})$ est 
donc de dimension 1 si $\ell_{1}\neq\ell_{2}$ et de dimension 2 si
$\ell_{1}=\ell_{2}$ puisque 
$\text{Hom}_{X}(I_{\ell_{1}},I_{\ell_{2}})=0$ 
si $\ell_{1}\neq\ell_{2}$ et 
$\text{Hom}_{X}(I_{\ell_{1}},I_{\ell_{2}})\simeq\mathbb{C}$ si 
$\ell_{1}=\ell_{2}$.\qed\\
\newline
\indent (4.4) Soient $N\ge 1$ un entier et $V$ un espace vectoriel 
complexe. Soient $Q$ le sch\'ema de Hilbert-Grothendieck 
param\'etrant les quotients 
$V\otimes\mathcal{O}_{X}(-N)\longrightarrow E$ sur $X$ de rang 2 et 
de classes de Chern $c_{1}(E)=0$, $c_{2}(E)=2$, $c_{3}(E)=0$ et 
$L$ la polarisation naturelle ([S]). Le groupe $G=PGL(V)$ agit sur $Q$ et 
une puissance convenable de $L$ est $G$-lin\'earis\'ee. Soit 
$Q_{c}^{ss}$ l'ouvert des points $PGL(V)$-semi-stables correspondants 
\`a des quotients sans torsion et $Q_{c}$ l'adh\'erence 
de $Q_{c}^{ss}$ dans $Q$. Lorsque l'entier $N$ et l'espace vectoriel $V$ 
sont convenablement choisis les propri\'et\'es suivantes sont 
satisfaites. L'application 
$V\otimes\mathcal{O}_{X}\longrightarrow E(N)$ induit un isomorphisme 
$V\simeq H^{0}(E(N))$ et 
$h^{i}(E(k))=0$ pour $k\ge N$ et $i\ge 1$ et pour tout quotient 
$E$ de $Q_{c}$. Le point $[E]\in Q_{c}$ est $PGL(V)$-semi-stable
si et seulement si le faisceau $E$ est semi-stable si 
et seulement si $[E]\in Q^{ss}_{c}$. Le stabilisateur de $[E]$ dans 
$GL(V)$ s'identifie au groupe des automorphismes du faisceau $E$ et 
l'espace des modules $M$ est alors le quotient de Mumford :
$$Q_{c}^{ss}//G$$
\textsc{Lemme 4.5}.$-$\textit{Sous les hypoth\`eses 
pr\'ec\'edentes, le sch\'ema $Q_{c}^{ss}$ est lisse.}\\
\newline
\textit{D\'emonstration}.$-$L'espace tangent \`a $Q^{ss}_{c}$ en 
un point $[E]$ est isomorphe \`a $\text{Hom}_{X}(F,E)$ o\`u $F$ est le noyau 
de l'application 
$V\otimes\mathcal{O}_{X}(-N)\longrightarrow E$. Le sch\'ema $Q_{c}^{ss}$ 
est lisse en ce point si $\text{Ext}^{1}_{X}(F,E)=0$. 
Consid\'erons la suite exacte :
$$\text{Ext}_{X}^{1}(V\otimes\mathcal{O}_{X}(-N),E)\longrightarrow
\text{Ext}_{X}^{1}(F,E)\longrightarrow\text{Ext}_{X}^{2}(E,E)$$
On en d\'eduit une inclusion
$\text{Ext}_{X}^{1}(F,E)\subset\text{Ext}_{X}^{2}(E,E)$ puisque 
$h^{1}(E(N))=0$. Il suffit donc de prouver 
$\text{Ext}^{2}_{X}(E,E)=0$.
Si $E$ est stable et localement 
libre alors le r\'esultat est 
d\'emontr\'e par [M-T] (lemme 2.7). Si $E$ est stable non localement libre 
alors l'annulation cherch\'ee est donn\'ee par le lemme 4.2. 
Si $E$ est strictement semi-stable alors $E$ 
est extension de $I_{\ell_{2}}$ par $I_{\ell_{1}}$ o\`u 
$\ell_{1}\subset X$ et $\ell_{2}\subset X$ sont deux droites. 
L'annulation cherch\'ee r\'esulte alors facilement du lemme 4.3.\qed\\
\newline
\textsc{Th\'eor\`eme 4.6}.$-$\textit{Soit 
$X\subset\mathbb{P}^{4}$ une hypersurface cubique lisse. 
L'espace des modules $M_{X}$ des faisceaux semi-stables de rang 2 sur 
$X$ de classes de chern $c_{1}=0$, $c_{2}=2$ et $c_{3}=0$ est lisse 
de dimension 5.}\\
\newline
\textit{D\'emonstration}.$-$Soient $x\in Q_{c}^{ss}$ et 
$E$ le faisceau correspondant. Soit $Q_{c}^{s}\subset Q_{c}$ 
l'ouvert des faisceaux stables et $M_{X}^{s}$ l'ouvert des 
classes d'isomorphismes de faisceaux stables. Le sch\'ema  
$Q_{c}^{s}$ est un fibr\'e principal sous $G$ 
au dessus de $M_{X}^{s}$ et $M_{X}^{s}$ est donc lisse (4.5).
Il nous reste \`a \'etudier $M_{X}$ en 
$E=I_{\ell_{1}}\oplus I_{\ell_{2}}$ o\`u $\ell_{1}\subset X$ et 
$\ell_{2}\subset X$ sont deux droites (3.5). L'orbite 
$O(x)$ du point $x$ sous $G$ est ferm\'ee. Son stabilisateur 
$G_{x}$ est un groupe r\'eductif et il existe un 
sous-sch\'ema affine $W\subset Q_{c}^{ss}$
passant par $x$ localement ferm\'e et 
stable sous l'action de $G_{x}$ tel que le morphisme  
$W//G_{x}\longrightarrow Q_{c}^{ss}//G$ 
soit \'etale ([L]). Soit $(W,x)$ le germe de $W$ en $x$ et soit $F$ 
la restriction \`a $X\times (W,x)$ du quotient tautologique sur 
$X\times Q$. Alors $((W,x),F)$ est un espace de d\'eformation 
verselles pour le faisceau $E$ ([O] prop. 1.2.3).
Le germe $W$ est donc lisse en $x$ (4.3) et puisque le 
morphisme $W//G_{x}\longrightarrow Q_{c}^{ss}//G$
est \'etale il suffit  
donc prouver que le quotient $W//G_{x}$ est lisse en $[x]$. Or il existe 
un morphisme $G_{x}$-lin\'eaire $(W,x)\longrightarrow (T_{x}W,0)$
\'etale en $x$ tel que le morphisme induit 
$W//G_{x}\longrightarrow T_{x}W//G_{x}$ soit 
\'etale en $[x]$ ([L]). Il suffit donc de prouver que le quotient 
$T_{x}W//G_{x}$ est lisse en $0$.\\
\indent Supposons les droites 
$\ell_{1}$ et $\ell_{2}$ distinctes. L'espace tangent 
$T_{x}W=\text{Ext}_{X}^{1}(E,E)=
\oplus_{i,j}\text{Ext}^{1}_{X}(I_{\ell_{i}},I_{\ell_{j}})$ 
est de 
dimension 6 (4.3) et $G_{x}=G_{m}\times G_{m}$ agit 
sur ledit espace par la formule ([O] lemme 1.4.16):
$$(t,t').(\sum_{i,j}e_{i,j})=e_{1,1}+\frac{t}{t'}e_{1,2}
+\frac{t'}{t}e_{2,1}+e_{2,2}$$
On v\'erifie facilement que le 
quotient $T_{x}W//G_{x}$ est isomorphe \`a l'espace affine 
$\mathbb{A}^{5}$ et en particulier lisse en $0$.\\
\indent Supposons les droites $\ell_{1}$ et $\ell_{2}$ confondues 
et notons $\ell$ cette droite. L'espace tangent 
$T_{x}W=\text{Ext}_{X}^{1}(E,E)$ est de dimension 8 (4.3) et 
$G_{x}=PGL(2)$. Posons 
$T=\text{Ext}^{1}_{X}(I_{\ell},I_{\ell})$ et soit $U$ un espace vectoriel 
de dimension 2. Le groupe $G_{x}$ agit sur $T_{x}W=T\otimes\text{End}(U)$ 
par conjugaison sur $\text{End}(U)$ ([O] lemme 1.4.16). Le quotient 
$T_{x}//G_{x}$ est isomorphe \`a l'espace affine $\mathbb{A}^{5}$ 
([La] III cas 2) et en particulier lisse en $0$.\qed\\
\newline
\textsc{Lemme 4.7}.$-$\textit{Soient $X\subset\mathbb{P}^{4}$ une cubique 
lisse et $M_{X}$ l'espace des modules des faisceaux de rang 2 semi-stables de 
classes de Chern $c_{1}=0$, $c_{2}=2$ et $c_{3}=0$.  
Le sous-sch\'ema localement ferm\'e de $M_{X}$ 
des faisceaux stables non localement libres est 
irr\'eductible de dimension 4 et le sous-sch\'ema ferm\'e de $M_{X}$ des 
faisceaux strictement semi-stables est \'egalement irr\'eductible 
de dimension 4.}\\
\newline
\textit{D\'emonstration}.$-$Soient $B$ la surface de Fano de $X$ 
et $Z\subset X\times B$ la vari\'et\'e d'incidence. La vari\'et\'e $Z$ est
lisse et irr\'educible de dimension 3 et le morphisme 
$Z\longrightarrow B$ induit par la seconde projection fait de $Z$ un 
fibr\'e en droites projectives localement trivial pour la topologie 
de Zariski ([T]). Notons $X_{Z}$ la 
vari\'et\'e obtenue en \'eclatant $Z$ dans le produit $X\times B$ et 
notons $p$ et $q$ les morphismes induits sur $X$ et $B$ respectivement. 
La fibre du morphisme $q$ au-dessus d'un point $[\ell]\in B$ s'identifie 
\`a $X_{\ell}$ (1.1). Soit $Q$ le fibr\'e de rang 3 sur $B$ dont 
la fibre au-dessus d'un point $[\ell]\in B$ est l'ensemble des 
\'equations des hyperplans 
de $\mathbb{P}^{4}$ contenant la droite $\ell$. Le morphisme 
surjectif naturel $q^{*}Q\twoheadrightarrow p^{*}\mathcal{O}_{X}(1)$ induit 
un morphisme propre et dominant 
$X_{Z}\overset{c}{\longrightarrow}\mathbb{P}_{B}(Q)$. Le morphisme 
$p\times c$ induit un plongement de $X_{Z}$ dans 
$X\times\mathbb{P}_{B}(Q)$ au-dessus de $\mathbb{P}_{B}(Q)$. L'ensemble 
des points de $\mathbb{P}_{B}(Q)$ est en bijection ensembliste 
avec l'ensemble des coniques trac\'ees sur $X$. Soit 
$U\subset \mathbb{P}_{B}(Q)$ l'ouvert des coniques lisses et soit  
$\pi$ la projection de $\mathbb{P}_{B}(Q)$ sur $B$. Soit 
$x\in U$ et $C_{x}=c^{-1}(x)\subset X$ la conique 
correspondante. La conique $C_{x}$ engendre un plan de 
$\mathbb{P}^{4}$ dont l'intersection r\'esiduelle avec $X$ est la droite 
$\pi(x)=[\ell_{x}]$. Soit 
$Y={(c^{-1}(U)\cap Z_{U})}_{\text{red}}
\subset c^{-1}(U)\subset U\times\mathbb{P}^{4}$ 
o\`u $Z_{U}=Z\times_{B}U$. Le morphisme induit $Y\longrightarrow U$ 
est alors fini de degr\'e 2. La fibre du 
morphisme pr\'ec\'edent au dessus d'un point $x\in U$ 
est ensemblistement l'intersection $C_{x}\cap\ell_{x}$. 
Supposons $Y$ irr\'eductible. Le sch\'ema 
$c^{-1}(U)\times_{U}Y$ poss\'ede une section au dessus de $Y$. 
Ladite section d\'etermine un morphisme quasi-fini 
$Y\longrightarrow M_{X}$ dont l'image est 
pr\'ecis\'ement le sous-sch\'ema localement ferm\'e des faisceaux 
stables non localement libres (3.4). Le lemme est donc d\'emontr\'e 
dans ce cas. Si $Y$ n'est pas irr\'eductible alors le morphisme 
$c$ poss\'ede une section au-dessus de l'ouvert $U$ et l'argument 
pr\'ec\'edent s'applique directement.\\
\indent Les faisceaux strictement semi-stables sont param\'etr\'es par les 
couples de droites de $X$ (3.5). Or on a un morphisme 
quasi-fini $B\times B\longrightarrow M_{X}$ dont l'image est le ferm\'e 
des faisceaux strictement semi-stables. La surface $B$ est 
irr\'eductible et ledit ferm\'e l'est donc \'egalement.\qed\\
\newline
\textsc{Th\'eor\`eme 4.8}.$-$\textit{Soient 
$X\subset\mathbb{P}^{4}$ une cubique 
lisse et $M_{X}$ l'espace des modules 
des faisceaux de rang 2 semi-stables de 
classes de Chern $c_{1}=0$, $c_{2}=2$ et $c_{3}=0$. Soit $B$ la surface de 
Fano de $X$. Alors $M_{X}$ est isomorphe 
\`a l'\'eclatement d'un translat\'e de la surface $-B$  
dans la jacobienne interm\'ediaire de $X$.}\\
\newline
\textit{D\'emonstration}.$-$Soit $U\subset M_{X}$
l'ouvert des fibr\'es vectoriels stables ; $M_{X}\setminus U$ est de 
dimension 4 (3.5 et 4.7). L'espace des modules 
$M_{X}$ est lisse de dimension 5 (4.6) et l'ouvert $U$ est donc dense.
La vari\'et\'e $M_{X}$ est donc irr\'eductible (2.6 et [I-M] cor. 5.1).\\
\indent L'espace des modules $M_{X}$ est le quotient 
de Mumford $Q_{c}^{ss}//G$ 
o\`u $Q_{c}^{ss}$ est lisse (4.5). Soit $\mathcal{E}$ une famille 
universelle sur $Q_{c}^{ss}\times X$. Fixons $t_{0}\in Q_{c}^{ss}$. 
L'application 
$$
\begin{array}{rcl}
Q_{c}^{ss} & \rightarrow & J(X) \\
t & \mapsto & c_{2}(\mathcal{E}_{t})-c_{2}(\mathcal{E}_{t_{0}})
\end{array}
$$
est alg\'ebrique (1.1) et \'equivariante sous 
l'action du groupe $G$. 
On en d\'eduit un morphisme que nous noterons 
encore $c_{2}$ de $M_{X}$ vers $J(X)$. Ce morphisme est birationnel 
([I-M] thm. 3.2) et induit un isomorphisme par restriction \`a 
l'ouvert des fibr\'es stables (2.6 et [M-T]).\\
\indent Les vari\'et\'es $M_{X}$ et $J(X)$ sont lisses et 
le lieu exceptionnel $D$ de $c_{2}$ est donc de 
codimension pure 1. Le diviseur $D$ a au plus deux composantes 
irr\'eductibles. La restriction de $c_{2}$ 
au diviseur des faisceaux strictement semi-stables 
est g\'en\'eriquement finie et le morphisme $c_{2}$ est 
donc un isomorphisme au point g\'en\'erique du diviseur consid\'er\'e.\\
\indent La grassmanienne des plans de $\mathbb{P}^{4}$ est rationnelle et 
les cubiques planes trac\'ees sur $X$ sont donc toutes rationnellement 
\'equivalentes. Soient $C_{0}$ et $C_{1}$ deux coniques trac\'ees sur 
$X$ et $\ell_{0}$ et $\ell_{1}$ les intersections r\'esiduelles 
respectives des plans des coniques avec $X$. On a donc 
$C_{1}=C_{0}+\ell_{0}-\ell_{1}$ dans $J(X)$. Si $E$ est un faisceau 
associ\'e \`a une conique $C\subset X$ (3.4) alors $c_{2}(E)=C$. 
On en d\'eduit que le diviseur adh\'erence des faisceaux  
stables non localement libres  
est contract\'e sur un translat\'e de $-B$ dans $J(X)$. Ce diviseur 
est irr\'eductible (4.7) et $M_{X}$ est donc isomorphe \`a 
l'\'eclatement d'un translat\'e de $-B$ dans $J(X)$ ([L] thm. 2).\qed\\ 
\newpage
\centerline{\textbf{R\'ef\'erences}}
$\ $
\newline
\noindent [A-C-G-H] E. Arbarello, M. Cornalba, P. Griffits, J. Harris, 
$\textit{Geometry of Algebraic curves I}$, Grundleheren der math. 
Wissenschaften, vol. 267, Springer Verlag (1985).\\
\newline
[H1] R. Hartshorne, $\textit{Stable vector bundles of rank 2 on }
\mathbb{P}^{3}$, Math. Ann. 238 (1978), 229-280.\\
\newline
[H2] R. Hartshorne, $\textit{Stable reflexives sheaves}$, Math. Ann. 254 
(1980), 121-176.\\
\newline
[I-M] A. Iliev, D. Markushevich, \emph{The Abel-Jacobi map 
for a cubic threefold and periods of Fano threefolds of degree 14}, 
Documenta Math. 5 (2000), 23-47.\\
\newline
[La] Y. Laszlo, $\textit{Local structure of the moduli space of vector 
bundles over curves}$, Comment. Math. Helvetici 71 (1996), 373-401.\\
\newline
[L] D. Luna, $\textit{Slices \'etales}$, M\'em. Soc. Math. France 33
(1973), 81-105.\\
\newline
[Lu] Z. Luo, $\textit{Factorization of birational morphisms of regular 
schemes}$, Math. Z. 212 (1993), 505-509.\\
\newline
[M-T] D. Markushevich, A.S. Tikhomirov, \emph{The Abel-Jacobi map 
of a moduli component of vectors bundles on the cubic threefold}, 
math. AG/9809140, to appear in J. Alg. Geom.\\
\newline
[M] M. Maruyama, $\textit{Boundedness of semi-stable sheaves of small 
ranks}$, Nagoya Math. J. 78 (1980), 65-94.\\
\newline
[Mum] D. Mumford, $\textit{Lectures on curves on an algebraic surface}$,
Annals of Math. Sudies 59, Princeton university Press (1966).\\
\newline
[Mu] J.P. Murre, $\textit{Some results on cubic threefolds}$ dans 
\emph{Classification of algebraic varieties and 
compact complex manifolds}, 
Springer-Verlag Lecture Notes 412 (1974), 140-164.\\
\newline
[O] K. O'Grady, $\textit{Desingularized moduli spaces of sheaves on a 
K3}$, J. Reine Angew Math. 512 (1999), 49-117.\\
\newline
[S] C. Simpson, $\textit{Moduli of representations of the fundamental 
group of a smooth variety}$, Publ. Math. IHES 79 (1994), 47-129.\\
\newline
[T] A.N. Tyurin, $\textit{The Fano surface of a nonsingular cubic in } 
\mathbb{P}^{4}$, Izv. Akad. Nauk. SSSR Ser. Mat. 34 (1970), 1200-1208.
\vspace{0.5cm}\\
\begin{center}
St\'ephane \textsc{Druel}\\
DMA-\'Ecole Normale Sup\'erieure\\
45 rue d'Ulm\\
75005 PARIS\\
e-mail: \texttt{druel@clipper.ens.fr}
\end{center}

\end{document}